\documentclass{amsart}

\usepackage{latexsym}
\usepackage{amssymb}
\usepackage{amsmath}
\usepackage{amsthm}
\usepackage{hyperref}
\usepackage{bigints}

\newcommand{\nablab}{\overline{\nabla}}

\newcommand{\lgra}{\longrightarrow}

\newcommand{\iid}{\mathrm{Id}\,}

\newcommand{\ddiv}{\mathrm{div}\,}

\newcommand{\lto}{\ensuremath{\longrightarrow}}

\newcommand{\Ss}{\mathbb{S}}

\newcommand{\R}{\mathbb{R}}

\newcommand{\function}[5]
{\begin{eqnarray*}\begin{array}{r@{}ccl}
 #1\;\colon\;  & #2 &\lto & #3 \\[.05cm]  
  & #4 &\longmapsto  & #5 
\end{array}\end{eqnarray*}
}

\newtheorem{example}{Exemples}[section]
\newtheorem{thm}{Theorem}[section]

\newtheorem{lemma}[thm]{Lemma}

\newtheorem{remark}[thm]{Remark}
\newtheorem{remarks}[thm]{Remarks}
\newtheorem{definition}[thm]{Definition}
\newtheorem{notation}[thm]{Notation}
\newtheorem{exabout:ample}[thm]{Example}

\newcommand{\beqt}{\begin{equation}}  \newcommand{\eeqt}{\end{equation}}
\newcommand{\bal}{\begin{align}}      \newcommand{\eal}{\end{align}}
\newcommand{\ba}{\begin{array}}      \newcommand{\ea}{\end{array}}
\newcommand{\bc}{\begin{center}}     \newcommand{\ec}{\end{center}}
\newcommand{\be}{\begin{enumerate}}  \newcommand{\ee}{\end{enumerate}}
\newcommand{\beq}{\begin{eqnarray}}  \newcommand{\eeq}{\end{eqnarray}}
\newcommand{\beQ}{\begin{eqnarray*}} \newcommand{\eeQ}{\end{eqnarray*}}
\newcommand{\bi}{\begin{itemize}}    \newcommand{\ei}{\end{itemize}}
\newcommand{\bt}{\begin{tabular}}    \newcommand{\et}{\end{tabular}}


\date{\today}

\begin{document}
\author{Julien Roth}
\address{Laboratoire d'Analyse et de Math\'ematiques Appliqu\'ees, UPEM-UPEC, CNRS, F-77454 Marne-la-Vall\'ee}
\email{julien.roth@u-pem.fr}

\title{On the stability of the Wulff shape}

\maketitle

\begin{abstract}
Given a positive function $F$ on $\mathbb S^n$ satisfying an appropriate convexity assumption, we consider hypersurfaces for which a linear combination of some higher order anisotropic curvatures is constant. We define the variational problem for which these hypersurfaces are critical points and we prove that, up to translations and homotheties, the Wulff shape is the only stable closed hypersurface of the Euclidean space for this problem.
\end{abstract}

\section{Introduction}
The stability of hypersurfaces for volume preserving variational problem has a long history since the first result for the stability of constant mean curvature in the Euclidean space by Barbosa and do Carmo \cite{BdC}. Many authors have been interested in stability problems in various contexts, like for other space forms and/or  higher order mean curvatures (see \cite{AdCC,AdCR,BdCE,BC,CW,Gr,LV,HL,HL3} for instance). In 2013, Vel\'asquez, de Sousa and de Lima \cite{VSL} defined the notion of $(r,s)$-stability which generalizes the classical notion of stability for the mean curvature or $r$-stability for higher order mean curvatures. They prove that the only compact $(r,s)$-stable hypersurfaces in the sphere or the hyperbolic space are geodesic hyperspheres. This result was recently extended for hypersurfaces of the Euclidean space by da Silva, de Lima and Vel\'asquez \cite{SLV}.\\ \\
On the other hand, during the last decade, an intensive interest has been brought to the study of hypersurfaces of Euclidean spaces in an anisotropic setting. Many of the classical characterizations of the geodesic hyperspheres have an analogue with the Wulff Shape as characteristic hypersurface, like anisotropic Hopf or Alexandrov-type theorems (see \cite{HL2,HL4,HLMG,KP,On,Pa}). In particular, in \cite{On}, Onat proved that a closed convex hypersurface of the Euclidean space with linearly related anisotropic mean curvatures $H^F_r=a_0H_0^F+a_1H^F+a_2H^F_2+\cdots+a_{r-1}H^F_{r-1}$ is the Wulff shape. The hypothesis that $X(M)$ is convex is crucial in Onat's result.  \\ \\
The aim of this short note is to prove an anisotropic analogue to \cite{SLV}. This extend to the anisiotropic $(r,s)$-stability the results of \cite{Pa}, \cite{HL4} and \cite{CdS} and give an other characterization of the Wulff shape as the only hypersurface (up to translations and homotheties) which have linearly related anisotropic mean curvatures, without assuming that $X(M)$ is convex. Namely, we prove the following
\begin{thm}\label{thm1}
Let $n,r,s$ three integers satisfying $0\leqslant r\leqslant s\leqslant n-2$. Let $F:\Ss^n\lgra\R_+$ be a smooth function satisfying the following convexity assumption  $$A_F=(\nabla dF+F\iid_{|T_x\Ss^n})_x>0,$$
for all $x\in\Ss^n$ and let $X:M^n\lgra\R^{n+1}$ be a closed hypersurface with positive anisotropic $(s+1)$-th mean curvature $H_{s+1}^F$. Assume that the quantity $\sum_{j=r}^s a_jb_jH^F_{j+1}$ is constant, where $a_r,\cdots,a_s$ are some nonnegative constants (with at least one non zero) and $b_j=(j+1)\binom{n}{j+1}$ for any $j\in\{r,\cdots,s\}$. \\
Then $X:M^n\lgra\R^{n+1}$ is $(r,s,F)$-stable if and only $X(M)$ is the Wulff shape $W_F$, up to translations and homotheties.
\end{thm}
The notion of $(r,s,F)$-stability will be defined in Section \ref{sec22}.
\section{Preliminaries}

\subsection{Anisotropic mean curvatures} 
Here, we recall the basics of anisotropic mean curvatures. These facts are classical, hence, we will not recall their proofs. First, let $F:\Ss^n\lgra\R_+$ be a smooth function satisfying the following convexity assumption
\begin{equation}\label{convexity}
A_F=(\nabla dF+F\iid_{|T_x\Ss^n})_x>0,
\end{equation}

Let $(M^n,g)$ be a compact, connected, oriented Riemannian manifold without boundary isometrically immersed into $\R^{n+1}$ by $X$ and denote by $\nu$ a normal unit vector field. Let $X^T=X-<X,\nu>\nu$ be its projection on the tangent bundle of $X(M)$. \\

The (real-valued) second fundamental form $B$ of the immersion is defined by
$$B(X,Y)=\left\langle \nablab_X\nu,Y\right\rangle ,$$
where $<\cdot,\cdot>$ and $\nablab$ are respectively the Riemannian metric and the Riemannian connection of $\R^{n+1}$. We also denote by $S$ the Weingarten operator, which is the $(1,1)$-tensor associated with $B$.\\ \indent
Let $S_F=A_F\circ S$. The operator $S_F$ is called the $F$-Weingarten operator and its eigenvalues $\kappa_1,\cdots,\kappa_n$ are the anisotropic principal curvatures. Now let us recall the definition of the anisotropic high order mean curvature $H_r^F$. First, we consider  an orthonormal frame $\{e_1,\cdots,e_n\}$ of $T_xM$. For all 
$k\in\{1,\cdots,n\}$, we set
$$\sigma_r=\left(\begin{array}{c}n\\r\end{array}\right)^{-1}\sum_{\begin{array}{c}1\leqslant i_1,\cdots,
    i_r\leqslant n\\1\leqslant j_1,\cdots,
    j_r\leqslant n\end{array}}\epsilon\left(\begin{array}{c}i_1\cdots
    i_r\\j_1\cdots
    j_r\end{array}\right)S^F_{i_1j_1}\cdots S^F_{i_rj_r},$$
    where the $S^F_{ij}$ are the coefficients of the $F$-Weingarten operator. The symbols $\epsilon\left(\begin{array}{c}i_1\cdots
    i_r\\j_1\cdots
    j_r\end{array}\right)$ are the usual permutation symbols which are
zero if the sets $\{i_1,\cdots,i_r\}$ and $\{j_1,\cdots,j_r\}$ are
    different or if there exist distinct $p$ and $q$ with $i_p=i_q$. For all other cases, $\epsilon\left(\begin{array}{c}i_1\cdots
    i_r\\j_1\cdots
    j_r\end{array}\right)$ is the signature of the permutation $\left(\begin{array}{c}i_1\cdots
    i_r\\j_1\cdots
    j_r\end{array}\right)$. 
    Then, the $r$-th anisotropic mean curvature of the immersion is defined by
    $$H_r^F=\left(\begin{array}{c}n\\r\end{array}\right)^{-1}\sigma_r.$$
We denote simply $H^F$ the anisotropic mean curvature $H_1^F$. Moreover, for convenience, we set $H^F_0=1$ et $H^F_{n+1}=0$ by convention.\\

For $r\in\{1,\cdots,n\}$, the symmetric $(1,1)$-tensor associated to $H_r^F$ is
$$P_r=\frac{1}{r!}\sum_{\begin{array}{c}1\leqslant i,i_1,\cdots,
    i_r\leqslant n\\1\leqslant j,j_1,\cdots,
    j_r\leqslant n\end{array}}\epsilon\left(\begin{array}{c}i,i_1\cdots
    i_r\\j,j_1\cdots
    j_r\end{array}\right)S^F_{i_1j_1}\cdots S^F_{i_rj_r}e_i^*\otimes e_j^*.$$
    We also define the following useful operators $T_r=P_rA_F$. Note that $T_r$ is symmetric. Moreover, we have these classical facts about the anisotropic mean curvatures (see \cite{HL} for instance).
    \begin{lemma}\label{lemPr} 
    For any $r$, we have
    \begin{enumerate}
    \item $P_r$ is divergence-free,
    \item $tr(P_r)=(n-r)\sigma_r$,
    \item $tr(P_rS_F)=(r+1)\sigma_{r+1}$,
    \item $tr(P_rS_F^2)=\sigma_1\sigma_{r+1}-\sigma_{r+2}$.
    \end{enumerate}    
    \end{lemma}
    
    \begin{lemma}\label{lemHF}
    Let  $r\in\{1,\cdot n-1\}$. If $H_{r+1}>0$, then For all $j\in\{1,\cdots,r\}$, 
    \begin{enumerate}
    \item $H^F_j>0$,
    \item $H^FH^F_{j+1}-H^F_{j+2}\geqslant 0$. Moreover, equality occurs at a point $p$ if and only if all the anisotropic principal curvatures at $p$ are equal. Hence, equality occurs everywhere if and only if $M$ is the Wulff shape, up to translations and homotheties.
    \end{enumerate}
     
    \end{lemma}
    Finally, we recall the anisotropic analogue of the classical Hsiung-Minkowski formulas \cite{Hs}. The proof can be found in \cite{HL} and uses in particular the fact that $P_r$ is divergence-free.
    \begin{lemma}\label{HS}
   Let  $r\in\{0,\cdots n-1\}$. Then, we have
   $$\int_M\left(F(\nu)H_r^F+H_{r+1}^F\langle X,\nu\rangle\right)dv_g.$$
    \end{lemma}
\subsection{The variational problem}\label{sec22}
In this section, we describe the stability problem that we will consider. For this we introduce the anisotropic $r$-area functionals
\beqt
\mathcal{A}_{r,F}=\left(\int_MF(\nu)\sigma_rdv_g\right), 
\eeqt
for $r\in\{0,\dots,n-1\}$ and where $dv_g$ denotes the Riemannian volume form on $M$.
Now we consider a variation of the immersion $X$. Precisely,  let $\varepsilon>0$ and 
\beqt 
\mathcal X\colon (-\varepsilon,\varepsilon)\times M\lgra\R^{n+1},
\eeqt
 such that for all $t\in (-\varepsilon,\varepsilon)$, $\mathcal X_t:=\mathcal X(t,\cdot)$ is an immersion of $M$ into $\R^{n+1}$ and $\mathcal X(0,\cdot)=X$. We denote by $\sigma_r(t)$ the corresponding curvature functions, by $\mathcal{A}_{r,F}(t)$ the $r$-area of  $\mathcal X_t$ and finally we set 
 \beqt f_t=\left\langle \frac{dX}{dt},\nu_t \right\rangle,\eeqt
  where $\nu_t$ is the unit normal to $M$ induced by $\mathcal X_t$. Finally, we denote by $g_t$ the induced metric on $M$.
  
Note that 
\beqt\label{firstvar}
\mathcal{A}_{r,F}'(t)=-b_{r+1}\int_M fH^F_{r+1}(t)dv_{g_{t}},\eeqt
where $b_{r+1}=(r+1)\binom{n}{r+1}$ cf.~\cite{HL4}.
We also consider the volume functional
\beqt
V(t)=\int_{[0,t)\times M}\mathcal X^*dv.
\eeqt
It is easy to see, cf. \cite[Lemma~2.1]{BdCE}, that $V$ satisfies
\beqt\label{Vprim}
 V'(t)=\int_Mf_tdv_{g_t}
\eeqt
and so $\mathcal X$ preserves the volume if and only if $\int_Mf_{t}dv_{g_t}=0$ for all $t$. Moreover, according to \cite[Lemma~2.2]{BdCE}, for any function $f_0:M\rightarrow\R$ such that $\int_Mf_0dv_g=0$, there exists a variation of $X$ preserving the volume and with  normal part given by $f_0$. \\
Now, let $r$ and $s$ two integers satisfying $0\leqslant r\leqslant s\leqslant n-2$ and $a_j$, $j=r,\cdots,s$ some nonnegative real numbers with at least one non zero. We consider the following anisotropic $(r,s)$-area functional $\mathcal{B}_{r,s,F}$ defined by 
\beqt\label{defBrsF}
\mathcal{B}_{r,s,F}=\sum_{j=r}^sa_j\mathcal{A}_{j,F}.
\eeqt
 This functional appears naturally when considering hypersurfaces with linearly related higher order anisotropic mean curvatures.  Indeed, we consider varaitions of $M$ that  preserve the balanced volume, the Jacobi functional associated with this anisotropic $(r,s)$-area is given by
\function{\mathcal{J}_{r,sF}}{(\varepsilon,\varepsilon)}{\R}{t}{\mathcal{B}_{r,s,F}(t)+\Lambda V(t),}
where $\Lambda$ is a constant to be determined. From \eqref{firstvar} and \eqref{Vprim}, we have imeediately
$$\mathcal{J}_{r,sF}'(t)=\int_M\left( -\sum_{j=r}^sa_jb_jH_{j+1}^F+\Lambda\right)f_tdv_{g_t}.$$
Hence, we have, like in the isotropic context, that $X$ satisfies $\displaystyle \sum_{j=r}^sa_jb_jH_{j+1}^F=constant$ if and only if $X$ it is a critical point of the functional $\mathcal{J}_{r,s,F}$, or equivalently, if and only if $X$ is a critical point of $\mathcal{B}_{r,s,F}$ for variations that preserve the balanced volume. Now, we give the definition of the anisotropic $(r,s)$-stability that we call $(r,s,F)$-stability.
\begin{definition}
Let $n,r,s$ three integers satistying $0\leqslant r\leqslant s\leqslant n-2$. Let $F:\Ss^n\lgra\R_+$ be a smooth function satisfying the following convexity assumption \eqref{convexity} and let $X:M^n\lgra\R^{n+1}$ be a closed hypersurface satisfying $$\sum_{j=r}^sa_jb_jH_{j+1}^F=constant.$$ Then, $X$ is said $(r,s,F)$-stable if $\mathcal{B}''_{r,s,F}(0)\geqslant 0$ for all volume-preserving variations of $X$.
\end{definition}
We consider the Jacobi operator $\mathcal{J}''_{r,s,F}(0)$ defined on the set $\mathcal F$ of smooth functions on $M$ with $\displaystyle\int_Mfdv_g=0$. 
From the definitions, we have clearly $\mathcal{B}''_{r,s,F}(0)=\mathcal{J}''_{r,s,F}(0)[f]$  where $f\in\mathcal F$ defines the variation $\mathcal X$. Therefore, the $(r,s,F)$-stability corresponds to the nonnegativity of the Jacobi operator.\\ \\
We finish this section by giving the second variation formula for this variational problem.

\begin{lemma}\label{secondvar}
For any variation of $mathcal X$ of $X$ preserving the balanced volume, the second variation formula of $\mathcal{B}_{r,s,F}$ at $t=0$ is given by
$$\mathcal{B}''_{r,s,F}(0)=\mathcal{J}''_{r,s,F}(0)[f]=-\sum_{j=r}^s(j+1)a_j\int_M\Big(L_jf+\langle T_j\circ d\nu,d\nu\rangle f\Big) fdv_g,$$
where $f\in\mathcal F$ is the normal part of the variation $\mathcal X$.
\end{lemma}
{\it Proof:} The proof comes directly from the second variation formula for each functional $\mathcal{A}_{j,F}$. Indeed, from \cite{HL4}, we have
$$\mathcal{A}''_{j,F}(0)=-(j+1)a_j\int_M\Big(L_jf+\langle T_j\circ d\nu,d\nu\rangle f\Big) fdv_g$$
Then, we have just to multiply by $a_j$ and sum from $r$ to $s$ to get the result. \hfill $\square$
\subsection{Some lemmas}
Now, we fix some notations and give some usefull lemmas. First, we define the following operators. For any $f\in\mathcal{C}^{\infty}$, we set
$$I_{j,F}[f]=L_jf+\langle T_j\circ d\nu,d\nu\rangle f$$
and 
$$\mathcal{R}_{r,s,F}=\sum_{j=r}^s(j+1)a_jI_{j,F}[f].$$
Obviously from this definition and Lemma \ref{secondvar}, we have
$$\mathcal{J}''_{r,s,F}(0)[f]=-\int_Mf\mathcal{R}_{r,s,F}[f]dv_g.$$
First, we recall this lemma due to He and Li \cite{HL4}. The proof of this lemma follows the idea of \cite{Ros} and uses Lemma \ref{lemPr}.
\begin{lemma}\label{lemI}
For any $j\in\{r,\cdots,s\}$, we have
\begin{enumerate}
\item $I_{j,F}[\langle X,\nu\rangle]=-\langle grad\,\sigma_{j+1},X^T\rangle-(j+1)\sigma_{j+1},$
\item $I_{j,F}[F(\nu)]=-\langle grad\,\sigma_{j+1},(grad_{\mathbb S^n}F)\circ\nu\rangle+\sigma_1\sigma_{j+1}-(j+2)\sigma_{j+2}$.
\end{enumerate}
\end{lemma}
Now, we have this last lemma about the symmetry of $\mathcal{R}_{r,s,F}$ w.r.t. the $L_2$-scalar product. Namely, we have
\begin{lemma}\label{lemR}
For any two smooth functions $f$ and $h$ over $M$, we have
$$\int_Mh\mathcal{R}_{r,s,F}[f]dv_g=\int_Mf\mathcal{R}_{r,s,F}[h]dv_g.$$
\end{lemma}
{\it Proof:} The proof is fairly standard. First, we compute 
\begin{eqnarray*}
\int_MhL_jfdv_g&=&\int_Mh \ddiv(T_j\nabla f)dv_g\\
&=&-\int_M \langle T_j\nabla f,\nabla h\rangle dv_g\\
&=&-\int_M \langle T_j\nabla h,\nabla f\rangle dv_g\\
&=&\int_Mf \ddiv(T_j\nabla h)dv_g\\
&=&\int_MfL_jhdv_g,
\end{eqnarray*}
where we have used the symmetry of $T_j$ and the divergence theorem. \\
Hence, from the defintion of $I_{j,F}$ and the above identity, we get immediately that
$$\int_MhI_{j,F}[f]dv_g=\int_MfI_{j,F}[h]dv_g.$$
Finally, multiplying by $a_jb_j$ and aking the sum over $j$ from $r$ to $s$, we get
$$\int_Mh\mathcal{R}_{r,s,F}[f]dv_g=\int_Mf\mathcal{R}_{r,s,F}[h]dv_g.$$
This concludes the proof.
\hfill $\square$\\ \\
Now, we have all the ingredients to prove the main result of this note.
\section{Proof of Theorem \ref{thm1}}
First, it is not difficult to see that the Wulff shape $W_F$ is $(r,s,F)$-stable. Indeed, the Wulff shape has all its mean curvatures $H_j^F$ constant and is $(j,F)$-stable, that is, $\mathcal{A}_{j,F}''(0)\geqslant 0$, for any $j\in\{0,\cdots,n-1\}$ (see \cite{HL4}). Therefore, $\displaystyle \sum_{j=r}^sa_jb_jH_{j+1}^F$ is clearly constant and since the constants $a_j$ are nonnegative, it is also clear that
$$\mathcal{B}_{r,s,F}''(0)=\sum_{j=r}^sa_j\mathcal{A}_{j,F}''(0)\geqslant0,$$ 
which means that the Wulff shape is $(r,s,F)$-stable.\\ \\
Conversely, suppose that $X:M^n\lgra\R^{n+1}$ is $(r,s,F)$-stable. By definition, we have $\mathcal{J}_{r,s,F}''(0)[f]\geqslant 0$ for any smooth function on $M$ satisfying $\int_M fdv_g=0$. We choose the particular test function $f$ defined by
$$f=\alpha F(\nu)+\beta\langle X,\nu\rangle,$$
with $ \alpha=\dfrac{\displaystyle\int_M\left(\sum_{j=r}^sa_jb_jF(\nu)H_j^F\right)dv_g}{\displaystyle\int_MF(\nu)}$ and $\beta=\displaystyle\sum_{j=r}^sa_jb_jH_{j+1}^F$. First, remark that $\beta$ is a constant by assumption. Moreover, using the anisotropic Hsiung-Minkowski formulas (Lemma \ref{HS}) we have
\begin{eqnarray*}
\int_Mfdv_g&=&\int_M\left(\alpha F(\nu)+\beta\langle X,\nu\rangle\right)dv_g\\
&=&\alpha\int_M F(\nu)dv_g+\int_M\left(\sum_{j=r}^sa_jb_jH_{j+1}^F\langle X,\nu\rangle\right)dv_g\\
&=&\alpha\int_M F(\nu)dv_g+\sum_{j=r}^sa_jb_j\int_MH_{j+1}^F\langle X,\nu\rangle dv_g\\
&=&\alpha\int_M F(\nu)dv_g-\sum_{j=r}^sa_jb_j\int_MF(\nu)H_{j}^F dv_g\\
&=&\alpha\int_M F(\nu)dv_g-\int_M\left(\sum_{j=r}^sa_jb_jF(\nu)H_{j}^F \right)dv_g\\
&=&\alpha\int_M F(\nu)dv_g-\alpha\int_M F(\nu)dv_g=0
\end{eqnarray*}
Hence, the integral of $f$ vanishes and $f$ is eligible as a test function. Hence, we have $\mathcal{J}_{r,s,F}''(0)[f]=-\int_Mf\mathcal{R}_{r,s,F}[f]dv_g\geqslant 0$. Now, let's compute $\mathcal{R}_{r,s,F}[f]$. We have

\begin{eqnarray*}
\mathcal{R}_{r,s,F}[f]&=&\sum_{j=r}^s(j+1)a_jI_{j,F}[f]\\
&=&\sum_{j=r}^s(j+1)a_jI_{j,F}[\alpha F(\nu)+\beta\langle X,\nu\rangle]\\
&=&\sum_{j=r}^s(j+1)a_j\Big( \alpha I_{j,F}[F(\nu)]+\beta I_{j,F}[\langle X,\nu\rangle]\Big).
\end{eqnarray*}
From Lemma \ref{lemI}, we have
\begin{eqnarray*}
\mathcal{R}_{r,s,F}[f]&=&\sum_{j=r}^s(j+1)a_j\Bigg[\alpha \Big(-\langle grad\,\sigma_{j+1},(grad_{\mathbb S^n}F)\circ\nu\rangle+\sigma_1\sigma_{j+1}-(j+2)\sigma_{j+2} \Big)\\
&&+\beta\Big( -\langle grad\,\sigma_{j+1},X^T\rangle-(j+1)\sigma_{j+1}\Big)\Bigg]
\end{eqnarray*}
Since $\sum_{j=r}^s (j+1)a_j\sigma_{j+1}=\sum_{j=r}^s a_jb_jH^F_{j+1}=constant,$ we get
\begin{eqnarray*}
\mathcal{R}_{r,s,F}[f]&=&\sum_{j=r}^s(j+1)a_j\Bigg[\alpha \Big(\sigma_1\sigma_{j+1}-(j+2)\sigma_{j+2}\Big) -\beta(j+1)\sigma_{j+1}\Bigg]\\
&=&\sum_{j=r}^sa_jb_j\Bigg[\alpha \Big(nH^FH^F_{j+1}-(n-j-1)H^F_{j+2}\Big) -\beta(j+1)H^F_{j+1}\Bigg]
\end{eqnarray*}

Moreover, we have
\begin{eqnarray}\label{int1}
\mathcal{J}_{r,s,F}''(0)[f]&=&-\int_Mf\mathcal{R}_{r,s,F}[f]dv_g  \nonumber \\
&=&-\int_M\Big(\alpha F(\nu)+\beta \langle X,\nu\rangle\Big)\mathcal{R}_{r,s,F}[f]dv_g \nonumber \\
&=&-\int_M \Big(\alpha F(\nu)\mathcal{R}_{r,s,F}[f]+\beta f\mathcal{R}_{r,s,F}[\langle X,\nu\rangle]\Big)dv_g .
\end{eqnarray}
where we have used Lemma \ref{lemR}. The first term in this intergal is
\begin{eqnarray}\label{int2}
\alpha F(\nu)\mathcal{R}_{r,s,F}[f]&=&\sum_{j=r}^sa_jb_jF(\nu)\Bigg[\alpha^2\Big(nH^FH^F_{j+1}-(n-j-1)H^F_{j+2}\Big)-\alpha\beta(j+1)H^F_{j+1}\Bigg] \nonumber\\
&=&\sum_{j=r}^sa_jb_jF(\nu)\Bigg[\alpha^2(n-j-1)\big(H^FH^F_{j+1}-H^F_{j+2}\big)+\alpha^2(j+1)H^FH^F_{j+1} \nonumber\\
&&-\alpha\beta(j+1)H^F_{j+1}\Bigg].
\end{eqnarray}
The second term in \eqref{int1} is, using Lemma \ref{lemI} and the fact that $\displaystyle \sum_{j=r}^sa_jb_jH^F_{j+1}$ is constant,
\begin{eqnarray}\label{int3}
\beta f\mathcal{R}_{r,s,F}[\langle X,\nu\rangle]&=&\beta f\sum_{j=r}^s(j+1)a_jI_{j,F}[\langle X,\nu\rangle] \nonumber\\
&=&\beta f\sum_{j=r}^s(j+1)a_j\Big( -\langle grad\,\sigma_{j+1},X^T\rangle-(j+1)\sigma_{j+1}\Big) \nonumber\\
&=&-\beta f\sum_{j=r}^sa_j(j+1)^2\sigma_{j+1} \nonumber\\
&=&-\beta f\sum_{j=r}^sa_jb_j(j+1)H^F_{j+1} \nonumber\\
&=&-\sum_{j=r}^sa_jb_j(j+1)\beta\Big(\alpha F(\nu)+\beta\langle X,\nu\rangle\Big)H^F_{j+1}.
\end{eqnarray}
Now, putting \eqref{int1} and \eqref{int2} in \eqref{int3}, we get
\begin{eqnarray}\label{int4}
\mathcal{J}_{r,s,F}''(0)[f]&=&-\sum_{j=r}^sa_jb_j(n-j-1)\alpha^2\int_MF(\nu)\left( H^FH^F_{j+1}-H^F_{j+2}\right)dv_g \nonumber\\
&&-\sum_{j=r}^sa_jb_j(j+1)\alpha^2\int_MF(\nu)H^FH^F_{j+1}dv_g\nonumber \\
&&+\sum_{j=r}^s2a_jb_j(j+1)\alpha\beta\int_MF(\nu)H^F_{j+1}dv_g \nonumber\\
&&+\sum_{j=r}^sa_jb_j(j+1)\beta^2\int_MH^F_{j+1}\langle X,\nu\rangle dv_g.
\end{eqnarray}
Using the anisotropic Hsiung-Minkowski formulas again, we have 
$$\sum_{j=r}^s2a_jb_j(j+1)\beta^2\int_MH^F_{j+1}\langle X,\nu\rangle dv_g =-\sum_{j=r}^s2a_jb_j(j+1)\beta^2\int_MF(\nu)H^F_j dv_g,$$
and therefore, \eqref{int4} becomes
 \begin{eqnarray}\label{int5}
\mathcal{J}_{r,s,F}''(0)[f]&=&-\sum_{j=r}^sa_jb_j(n-j-1)\alpha^2\int_MF(\nu)\left( H^FH^F_{j+1}-H^F_{j+2}\right)dv_g \nonumber\\
&&-\sum_{j=r}^sa_jb_j(j+1)\int_MF(\nu)\Big( H^FH^F_{j+1}\alpha^2-2H^F_{j+1}\alpha\beta+H^F_j\beta^2\Big)dv_g
\end{eqnarray}
Now, at a point $x$ on $M$, we consider the following second order polynomial
$$P_{j,F,x}(z)=F(\nu)\Big(H^FH^F_{j+1}z^2-2H^F_{j+1}\beta z+H^F_j\beta^2\Big).$$
The discriminant of $P_{j,F,x}$ is 
$$\Delta=4\beta^2F(\nu)^2\Big((H^F_{j_1})^2-H^FH^F_{j+1}H^F_j\Big)=4\beta^2F(\nu)^2H^F_{j+1}(H^F_{J+1}-H^FH^F_j).$$
Since, by assumption, $H_{s+1}>0$, by Lemma \ref{lemHF}, we have that $H^F_{j+1}>0$ and $H^F_{J+1}-H^FH^F_j\geqslant0$. Hence, $\Delta$ is nonnegative and since the term of degree $2$ is $F(\nu)H^FH^F_{j+1}>0$, then $P_{j,F,x}(z)\geqslant 0$ for any $z\in\R$. In particular, for $z=\alpha$, we obtain
$$H^FH^F_{j+1}\alpha^2-2H^F_{j+1}\alpha\beta+H^F_j\beta^2\geqslant0.$$
Reporting this in \eqref{int5}, we get
$$\mathcal{J}_{r,s,F}''(0)[f]\leqslant-\sum_{j=r}^sa_jb_j(n-j-1)\alpha^2\int_MF(\nu)\left( H^FH^F_{j+1}-H^F_{j+2}\right)dv_g.$$
Finally, since $F(\nu)>0$ and $H^FH^F_{j+1}-H^F_{j+2}\geqslant0$ by Lemma \ref{lemHF}, we get that
$$\mathcal{J}_{r,s,F}''(0)[f]\leqslant0.$$
Since, by the $(r,s,F)$-stability assumption, we have $\mathcal{J}_{r,s,F}''(0)[f]\geqslant0$, we deduce that $$\mathcal{J}_{r,s,F}''(0)[f]=0.$$
This means that each term in the right-hand side of \eqref{int5} vanishes. In particular, we have
$$\sum_{j=r}^sa_jb_j(n-j-1)\alpha^2\int_MF(\nu)\left( H^FH^F_{j+1}-H^F_{j+2}\right)dv_g=0,$$
and so for each $j\in\{r,\cdots,s\}$, 
$$\int_MF(\nu)\left( H^FH^F_{j+1}-H^F_{j+2}\right)dv_g=0.$$
Since $F(\nu)>0$ and  $H^FH^F_{j+1}-H^F_{j+2}\geqslant0$, we get that at any point $x$ of $M$,
$$H^FH^F_{j+1}-H^F_{j+2}=0.$$
Thus, by Lemma \ref{lemHF}, we conclude that $X(M)$ is the Wulff shape, up to translations and homotheties. This concludes the proof of Theorem \ref{thm1}. \hfill$\square$


\begin{thebibliography}{00}

\bibitem{AdCC} H. Alencar, M. do Carmo \& A.G. Colares, \emph{Stable hypersurfaces with constant scalar curvature}, Math. Z. {\bf 213} (1993), 117-131.
\bibitem{AdCR} H. Alencar, M.P. Do Carmo \& H. Rosenberg, \emph{On the first eigenvalue of Linearized operator of the r-th mean curvature of a hypersurface}, Ann. Glob. Anal. Geom., {\bf 11} (1993), 387-395.
\bibitem{BdC} J.L.M. Barbosa \& M. do Carmo, \emph{Stability of hypersurfaces with constant mean curvature}, Math. Z. {\bf 185} (1984), 339-353.
\bibitem{BdCE} J.L.M. Barbosa, M.P. Do Carmo \& J. Eschenburg, \emph{Stability of hypersurfaces with constant mean curvature in Riemannian manifolds}, Math. Z. {\bf 197} (1988), 123-138.
\bibitem{BC} J.L.M. Barbosa \& A.G. Colares, \emph{Stability of hypersurfaces with constant $r$-mean curvature}, Ann. Glob. Anal. Geom. {\bf 15} (1997), 277-297.
\bibitem{CdS} A.G. Colares \& J.F. da Silva, \emph{ Stable hypersurfaces as minima of the intergral of an anisotropic mean curvature preserving a linear combination of area and volume}, Math. Z. {\bf 275} (2013), 595-623.
\bibitem{CW} H. Chen \& X. Wang, \emph{stability and eigenvalue estimates of linear Weingarten hypersurfaces in a sphere}, J. Math. Anal. Appl. {\bf 397 (2)} (2013), 653-670.
\bibitem{Gr} J.F. Grosjean, \emph{Extrinsic upper bounds for the first eigenvalue of elliptic operators}, Hokkaido Math. J. {\bf 33} (2004), no. 2, 219-239.
\bibitem{LV} H.F. de Lima\& M.A.L. Vel\'asquez, \emph{A new characterization of $r$-stable hypersurfaces in space forms}, Arch. Math, (Brno) {\bf 47} (2011), 119-131.
\bibitem{HL} Y. He \& H. Li, \emph{Integral formula of Minkowski type and new characterization of the Wulff Shape}, Acta Math. Sin. {\bf 24} (2008), no. 3, 697-704.
\bibitem{HL2} Y. He \& H. Li, \emph{A new variational characterization of the Wulff shape}, Diff. Geom. Appl., {\bf 26} (2008), no. 4, 377-390.
\bibitem{HL3} Y. He \& H. Li \emph{Stability of area-preserving variations in space forms}, Ann. Glob. Anal. Geom. {\bf 34} (2008), 55-68.
\bibitem{HL4} Y. He \& H. Li, \emph{Stability of hypersurfaces with constant $(r+1)$-th anisotropic mean curvature}, Illinois J. Math. {\bf 52 (4)} (2008), 1301-1314.
\bibitem{HLMG} Y. He, H. Li, H. Ma \& J.Ge, \emph{Compact embedded hypersurfaces with constant higher order anisotropic mean curvature}, Indiana Univ. Math. J. {\bf 58} (20058), 853-868.
\bibitem{Hs} C.C. Hsiung, \emph{Some integral formulae for closed hypersurfaces}, Math. Scand {\bf 2} (1954), 286-294.
\bibitem{KP} M. Koiso and B. Palmer, \emph{Geometry and stability of surfaces with constant anisotropic mean curvature}, Indiana Univ. Math. J. {\bf 54} (2005), 1817-1852. 
\bibitem{On} L. Onat, \emph{Some characterizations of the Wulff shape}, C. R. Math. Acad. Sci. Paris {\bf 348} (2010), no. 17-18, 997Ð1000.
\bibitem{Pa} B. Palmer, \emph{Stability of the Wulff shape}, Proc. Amer. Math. Soc. {\bf 126} (1998), no. 2, 3661-3667.
\bibitem{Ros} H. Rosenberg, \emph{Hypersurfaces of constant curvature in space forms}, Bull. Sc. Math. {\bf 117} (1993), 221-239.
\bibitem{SLV} J.F. da Siva, H.F. de Lima \& M.A.L. Vel\'asquez, \emph{The stability of hypersurfaces revisited}, Monatsh. Math (in press) DOI 10.1007/s00605-015-0776-x.
\bibitem{VSL} M.A.L. Vel\'asquez, A.F. de Sousa \& H.F. de Lima, \emph{On the stability of hypersurfaces in space forms}, J. Math. Anal. Appl. {\bf 406} (2013), 134-146.








\end{thebibliography}
\end{document}